\documentclass[12pt,reqno]{amsart}

\usepackage{amssymb}
\usepackage{amsmath}
\usepackage{amsthm}
\usepackage{mathrsfs}

\numberwithin{equation}{section}

\theoremstyle{plain}
\newtheorem{thm}{Theorem}[section]
\newtheorem{prop}[thm]{Proposition}
\newtheorem{lem}[thm]{Lemma}

\theoremstyle{definition}
\newtheorem{df}{Definition}[section]

\newtheorem{rem}{Remark}[section]

\DeclareMathOperator{\al}{a.e.}
\DeclareMathOperator{\as}{as}
\DeclareMathOperator{\an}{and}
\DeclareMathOperator{\re}{Re}
\DeclareMathOperator{\kr}{Ker}

\def\a{\alpha}
\def\b{\beta}
\def\bd{D_t}
\def\cd{\d_t}
\def\d{\partial}
\def\de{\delta}
\def\dis{\displaystyle}
\def\e{\varepsilon}
\def\g{\gamma}
\def\l{\lambda}
\def\nn{\nonumber}

\def\ov{\overline}

\def\ds{d\sigma}
\def\th{\theta}
\def\tm{\times}
\def\v{\varphi_n}
\def\vk{\varphi_{kj}}
\def\vlj{\varphi_{\ell j}}
\def\wt{\widetilde}

\def\C{\mathbb{C}}
\def\D{\mathcal{D}}

\def\G{\Gamma}
\def\L{\mathscr{L}}
\def\Lc{\mathcal{L}}
\def\Lm{\Lambda}
\def\N{\mathbb{N}}
\def\O{\Omega}
\def\Ov{\overline{\O}}
\def\R{\mathbb{R}}

\def\sm{\sum_{n=1}^{\infty}}
\def\smN{\sum_{n=1}^N}
\def\smMN{\sum_{n=M}^N}

\def\smij{\sum_{i,j=1}^d}
\def\smk{\sum_{k=1}^{\infty}}
\def\smj{\sum_{j=1}^{m_k}}
\def\smjl{\sum_{j=1}^{m_\ell}}

\def\({\left(}
\def\){\right)}
\def\<{\left\langle}
\def\>{\right\rangle}

\def\li{\varliminf}

\allowdisplaybreaks

\begin{document}

\title[approximate controllability]
{Approximate controllability for fractional diffusion equations by 
Dirichlet boundary control}

\author{Kenichi Fujishiro}



\thanks{Department of Mathematical Sciences, University of Tokyo, 
Komaba, Meguro, Tokyo 153, Japan}

\email{kenichi@ms.u-tokyo.ac.jp}

\date{}

\begin{abstract}
  In this paper, we consider the approximate controllability of partial 
differential equations with time derivatives of non-integer order via boundary 
control.
  First we show the unique existence and regularity of the solution by using 
the eigenfunction expansion.
  Next we also study the dual system and show the unique continuation property.
  Finally we apply it to prove our main result.
\end{abstract}

\maketitle

\section{Introduction}\label{sec:int}

  Let $\O$ be a bounded domain of $\R^d$ with $C^2$ boundary $\G=\d\O$.
  We consider the following initial value/boundary value problem of partial 
differential equation:
\begin{equation}\label{eq:ori}
\begin{cases}
	\cd^{\a}u+\L u=0	&\mbox{in}\quad\O\tm(0,T),\\
	u=g			&\mbox{on}\quad\G\tm(0,T),\\
	u(\cdot,0)=u_0		&\mbox{in}\quad\O.
\end{cases}
\end{equation}
  In \eqref{eq:ori}, $u=u(x,t)$ is the state to be controlled and $g=g(x,t)$ 
is the control which is localized on a subboundary $\G_0$ of $\G$.
  We will act by $g$ to drive the initial state $u_0=u_0(x)$ to some target 
function $u_1=u_1(x)$.
  Here $\L$ denotes a symmetric and uniformly elliptic operator, which is 
specified later and $T>0$ is a fixed value.
  The Caputo fractional derivative $\cd^{\a}$ is defined by
\begin{equation}\label{eq:ca}
	\cd^{\a}h(t):=
\begin{cases}
	\dfrac{d^n h}{d t^n}(t),	&\a=n\in\N, \\
	\dis\frac{1}{\G(n-\a)}\int_0^t(t-\tau)^{n-\a-1}
		\frac{d^n h}{d\tau^n}(\tau)d\tau,	&n-1<\a<n,\ n\in\N,
\end{cases}
\end{equation}
  for $\a>0$ (see \cite{ks}, \cite{pd}).
  If $\a=1$, then equation \eqref{eq:ori} is a classical diffusion equation.
  Equation \eqref{eq:ori} with $0<\a<1$ is called a {\it fractional diffusion 
equation} and regarded as a model of anomalous diffusion in heterogeneous 
media.
  In the present paper, we consider the case of $0<\a<1$.

  According to Adams and Gelhar \cite{ag}, the field data in a highly 
heterogeneous aquifer cannot be described well by classical advection 
diffusion equations.
  Hatano and Hatano \cite{ha} applied the continuous-time random walk (CTRW) 
as a microscopic model of the diffusion of ions in heterogeneous media.
  From the CTRW model, Metzler and Klafter \cite{mk} derived equation 
\eqref{eq:ori} with $0<\a<1$ as a macroscopic model.
  Concerning the mathematical analysis of fractional differential equations, 
there are many works.
  For the study of ordinary differential equations with fractional orders, see 
\cite{ks}, \cite{pd} and \cite{sm} for example.
  As for partial differential equations with time fractional derivatives, we 
can refer to Gejji and Jafari \cite{gj}, Agarwal \cite{aw} and Luchko 
\cite{lu} for example.

  The aim of this article is to study the boundary control problem for 
fractional diffusion equations.
  We say that equation \eqref{eq:ori} is {\it approximately controllable} if 
for any $u_1\in L^2(\O)$ and $\e>0$, there exists a control $g$ such that the 
solution $u$ of \eqref{eq:ori} satisfies
\begin{equation}\label{eq:e>0}
	\|u(\cdot,T)-u_1\|_{L^2(\O)}\le\e.
\end{equation}
  We can refer to \cite{co}, \cite{mz} and \cite{ru} for the general theory of 
control problems for partial differential equations.
  These surveys deal with controllability of equations with integer order and 
the relations with other concepts---observability, stabilizability, pole 
assignability, etc.
  There are various works about control problems for equations with integer 
orders.
  In particular, for the boundary control of heat equations, see MacCamy, 
Mizel and Seidman \cite{ma}, Schmidt and Weck \cite{sw}, 
and the references therein.
  However, to the author's best knowledge, there are few works on the 
control problems for fractional diffusion equations, especially on the 
boundary control problems.

\bigskip
  The remainder of this paper is composed of five sections and an appendix.
  In Section \ref{sec:main}, we state the main result.
  In Section \ref{sec:reg}, we give a representation of the solution by 
Fourier's method and discuss its fundamental properties.
  In Section \ref{sec:weak}, we define the weak solution by using the 
representation obtained in Section \ref{sec:reg}.
  In Section \ref{sec:dual}, we study the dual system of \eqref{eq:ori} and 
prove the unique continuation property, which plays an essential role in the 
proof of 
our main result.
  In Section \ref{sec:proof}, we complete the proof of the main result.
  In the appendix, we discuss the related boundary value problem for an 
elliptic equation.

\section{Main result}\label{sec:main}

  In this section, we prepare the notations and state our main results.
  We denote by $L^p(\O)$, $1\le p\le\infty$, a usual $L^p$-space.
  In particular, $L^2(\O)$ denotes the $L^2$-space equipped with the scalar 
product $(\cdot,\cdot)$.
  As for the inner product in $L^2(\G)$, we denote it by $\<\cdot,\cdot\>$.
  Moreover $H^l(\O)$ and $H^m_0(\O)$, $l,m\in\N$, are the Sobolev spaces (see 
Adams \cite{ar} for example).
  In equation \eqref{eq:ori}, let the differential operator $\L$ be given by
\[
	\L u(x)
	=-\smij\frac{\d}{\d x_i}\(a_{ij}(x)\frac{\d u}{\d x_j}(x)\)+c(x)u(x),
	\quad x\in\O,
\]
  where the coefficients satisfy the following:
\begin{flalign*}
	&a_{ij}=a_{ji},\ a_{ij}\in C^1(\Ov),\quad 1\le i,j\le d,\qquad
	\smij a_{ij}(x)\xi_i\xi_j\ge\mu|\xi|^2,\quad x\in\Ov,\ \xi\in\R^d, \\
	&c\in C(\Ov),\qquad c(x)\ge0,\quad x\in\Ov,
\end{flalign*}
  where $\mu>0$ is constant.
  We define the operator $L:L^2(\O)\to L^2(\O)$ as $\L$ equipped with 
the homogeneous Dirichlet boundary condition;
\begin{align*}
	&\D(L):=H^1_0(\O)\cap H^2(\O), \\
	&L u:=\L u,\quad u\in\D(L).
\end{align*}
  Since $L$ is a symmetric and uniformly elliptic operator, the spectrum of 
$L$ is entirely composed of countable number of eigenvalues and we can set 
with finite multiplicities:
\[
	0<\l_1\le\l_2\le\cdots\le\l_n\le\cdots.
\]
  By $\v\in H^2(\O)\cap H^1_0(\O)$, we denote the orthonormal eigenfunction 
corresponding to $\l_n$:
\[
	L\v=\l_n\v,\quad n=1,2,\cdots.
\]
  Then the sequence $\{\v\}_{n\in\N}$ is an orthonormal basis in $L^2(\O)$.
  We can represent the fractional power of $L$ as follows;
\begin{align*}
	&\D(L^{\th})=
	\left\{u\in L^2(\O);\ \sm\l_n^{2\th}|(u,\v)|^2<\infty\right\}, \\
	&L^{\th}u=\sm\l_n^{\th}(u,\v)\v,	\quad u\in\D(L^{\th}),
\end{align*}
  where $\th>0$.
  Then $\D(L^{\th})$ is a Hilbert space equipped with the norm 
$\|\cdot\|_{\D(L^{\th})}$ defined by
\[
	\|u\|_{\D(L^{\th})}
:=	\|L^{\th}u\|_{L^2(\O)}
=	\(\sm\l_n^{2\th}|(u,\v)|^2\)^{1/2},	\quad u\in\D(L^{\th}).
\]
  The domain $\D(L^{\th})$ with $0\le\th\le1$, $\th\ne1/4$, is expressed by 
using the Sobolev spaces with norm equivalence;
\begin{align*}
	&\D(L^{\th})=
\begin{cases}
	H^{2\th}(\O), &0\le\th<1/4, \\
	H^{2\th}_D(\O), &1/4<\th\le1,
\end{cases} \\
	&C^{-1}\|u\|_{H^{2\th}}
\le	\|u\|_{\D(L^{\th})}
\le	C\|u\|_{H^{2\th}}, \quad u\in\D(L^{\th}),
\end{align*}
  where $H^s_D(\O):=\{u\in H^s(\O)\ |\ \g_0 u=0\}$ and the operator 
$\g_0:H^s(\O)\to H^{s-1/2}(\G)$ maps a function $u$ to its restriction 
$u|_{\G}$ to the boundary $\G$ for $s>1/2$.
  For the details of $\D(L^{\th})$ and the Sobolev spaces with fractional 
powers, see Fujiwara \cite{fu} and Yagi \cite{yg} for example.
  The operator $\d_{\nu_L}:H^s(\O)\to H^{s-3/2}(\G)$, $s>3/2$, is defined as
\[
	\d_{\nu_L}u(x)
	=\frac{\d u}{\d\nu_L}(x)
	=\smij a_{ij}(x)\frac{\d u}{\d x_i}(x)\nu_j(x),
\]
  where $\nu(x)=(\nu_1(x),\dots,\nu_d(x))$ is the outward unit normal vector 
to $\G$ at $x$.
  In particular, $\d_{\nu_L}\v$ belongs to $L^2(\G)$ since $\v\in H^2(\O)$.
  We define the Mittag-Leffler function by
\[
	E_{\a,\beta}(z):=\sum_{k=0}^{\infty}\frac{z^k}{\G(\a k+\beta)},
	\quad z\in\C,
\]
  where $\a>0$ and $\beta\in\R$ are arbitrary constants. 
  We can directly verify that $E_{\a,\beta}(z)$ is an entire function of 
$z\in\C$.
  Henceforth $C$ denotes the positive generic constant which is independent of 
$g$, but may depend on $\a$ and the coefficients of the operator $L$.

  According to Theorem 2.1 in \cite{sk} and Proposition \ref{prop:reg1} in the 
next section, for any $u_0\in L^2(\O)$ and $g\in C_0^{\infty}(\G_0\tm(0,T))$, 
equation \eqref{eq:ori} admits a unique solution $u\in C([0,T];L^2(\O))$ with 
the representation as;
\begin{align}\label{eq:sol}
	u(x,t)=
&	\sm(u_0,\v)E_{\a,1}(-\l_n t^{\a})\v(x) \nn \\
&\quad	-\sm\(\int_0^t\<g(\cdot,t-\tau),\d_{\nu_L}\v\>
		\tau^{\a-1}E_{\a,\a}(-\l_n\tau^{\a})d\tau\)\v(x).
\end{align}
  In particular, the value $u(\cdot,T)$ at time $t=T$ makes sense in $L^2(\O)$ 
and consequently we can discuss problems such as whether \eqref{eq:e>0} is 
possible or not.

  Now we are ready to state one of our main results;

\bigskip
\begin{thm}\label{thm:ac}
  Let $0<\a<1$ and $\G_0$ be an open set in $\G$.
  Then equation \eqref{eq:ori} is approximately controllable for arbitrarily 
given $T>0$.
  That is,
\begin{equation}\label{eq:ac}
	\ov{\{u(\cdot,T);\ g\in C_0^{\infty}(\G_0\tm(0,T))\}}=L^2(\O),
\end{equation}
  where $u$ is the solution to \eqref{eq:ori} and the closure on the 
left-hand side is taken in $L^2(\O)$.
\end{thm}

\bigskip
  In order to prove this theorem, we also need to consider the {\it dual 
system} for \eqref{eq:ori}, which is a usual strategy for partial differential 
equations of integer order (see Section 8 in \cite{ru} or Chapters 2 and 3 in 
\cite{mz} for example).
  The dual system for \eqref{eq:ori} corresponds to the following initial 
value/boundary value problem with a different type of fractional derivative;
\begin{equation}\label{eq:dual}
\begin{cases}
	\bd^{\a}v+\L v=0		&\mbox{in}\quad\O\tm(0,T),\\
	v=0				&\mbox{on}\quad\G\tm(0,T),\\
	I_{T-}^{1-\a}v(\cdot,T)=v_0	&\mbox{in}\quad\O.
\end{cases}
\end{equation}
  Here $\bd^{\a}$ denotes the backward Riemann-Liouville derivative and is 
defined by
\begin{equation}\label{eq:rl}
	\bd^{\a}h(t)=
\begin{cases}
	\dis\(-\frac{d}{dt}\)^n h(t),	& \a=n\in\N, \\
	\dis\frac{1}{\G(n-\a)}\(-\frac{d}{dt}\)^n\int_t^T 
		(\tau-t)^{n-\a-1}h(\tau)d\tau,\ 	& n-1<\a<n,\ n\in\N,
\end{cases}
\end{equation}
  for $\a>0$ (see \cite{pd}).
  Moreover $I_{T-}^{\a}$ is the backward integral operator, which is defined by
\[
	I_{T-}^{\a}h(t):=\frac{1}{\G(\a)}\int_t^T(\tau-t)^{\a-1}h(\tau)d\tau
\]
  for $\a>0$ (see Section \ref{sec:reg} for details).
  In particular, for $0<\a<1$, we have
\begin{equation}\label{eq:frac}
	 \bd^{\a}h(t)=-\frac{d}{dt}I_{T-}^{1-\a}h(t).
\end{equation}
  We also note that the third equation in \eqref{eq:dual} means that
\[
	I_{T-}^{1-\a}v(x,T)
	:=\lim_{t\to T}\frac{1}{\G(1-\a)}\int_t^T(\tau-t)^{-\a}v(x,\tau)d\tau
	=v_0(x),\quad 0<\a<1.
\]

  In Section \ref{sec:dual}, we will study problem \eqref{eq:dual}.
  In Section \ref{sec:proof}, we will see that the unique continuation 
property for \eqref{eq:dual} is equivalent to the approximate controllability 
for \eqref{eq:ori} stated in Theorem \ref{thm:ac}.
  Moreover, by the variational approach, we can construct the control $g$ 
and show that it is also finite-approximately controllable.

  Let $E$ be a finite dimensional subspace of $L^2(\O)$ and fix $\e>0$ and 
$u_1\in L^2(\O)$ arbitrarily.
  We introduce the functional $J_{\e}$ on $L^2(\O)$ defined by
\begin{align}\label{eq:Je}
	J_{\e}(v_0):=
	\frac{1}{2}\int_0^T\int_{\G_0}(T-t)^2|\d_{\nu_L}v|^2\ds_xdt
	+&\e\|(I-\pi_E)v_0\|_{L^2(\O)} \nn \\
	+(v_0,u_1)&-(I_{T-}^{1-\a}v(\cdot,0),u_0),\quad v_0\in L^2(\O),
\end{align}
  where $v$ is the solution of \eqref{eq:dual} and 
$\pi_E$ denotes the orthogonal projection to $E$.
  By Proposition 4.1 in \cite{fy}, for any $v_0\in L^2(\O)$ equation 
\eqref{eq:dual} posesses a unique solution $v$ with 
$I_{T-}^{1-\a}v\in C([0,T];L^2(\O))$.
  Moreover $v$ is represented by
\[
	v(x,t)=\sm(v_0,\v)(T-t)^{\a-1}E_{\a,\a}(-\l_n (T-t)^{\a})\v(x).
\]
  Therefore $I_{T-}^{1-\a}v(\cdot,0)$ belongs to $L^2(\O)$ and 
$(T-t)^2|\d_{\nu_L}v|^2$ is integrable in $\G\tm(0,T)$.
  Thus the functional $J_{\e}$ is well-defined.

  Then we obtain the following result;

\bigskip
\begin{thm}\label{thm:ch}
  The functional $J_{\e}$ defined in \eqref{eq:Je} has a unique minimizer 
$\ov{v}_0\in L^2(\O)$.
  Moreover, let $\ov{v}$ be the solution of \eqref{eq:dual} with 
$v_0=\ov{v}_0$, then the solution $u$ of \eqref{eq:ori} with 
$g=(T-t)^2\d_{\nu_L}\ov{v}$ satisfies
\[
	\|u(\cdot,T)-u_1\|_{L^2(\O)}\le\e
	\quad\an\quad
	\pi_E(u(\cdot,T))=\pi_E(u_1).
\]
\end{thm}

\bigskip
  In the above theorem, we take $g=(T-t)^2\d_{\nu_L}\ov{v}$ as the control.
  However, in order to do this, we have to verify that \eqref{eq:ori} 
has a solution in $C([0,T];L^2(\O))$ for non-smooth $g$.
  In Section \ref{sec:weak}, therefore, we will define the {\it weak solution} 
of \eqref{eq:ori} for $g\in L^p(0,T;L^2(\O))$ with large $p$ and study its 
regularity.

  As for the variational method introduced here, we can refer to Lions 
\cite{li} and Zuazua \cite{zu}.

  We finally note that in Theorems \ref{thm:ac} and \ref{thm:ch}, we may 
assume $u_0=0$ without loss of generality.
  Indeed, consider the following two problems
\begin{equation}\label{eq:sim}
\begin{cases}
	\cd^{\a}u+\L u=0	&\mbox{in}\quad\O\tm(0,T),\\
	u=g			&\mbox{on}\quad\G\tm(0,T),\\
	u(\cdot,0)=0		&\mbox{in}\quad\O
\end{cases}
\end{equation}
  and
\begin{equation}\label{eq:homo}
\begin{cases}
	\cd^{\a}\wt{u}+\L\wt{u}=0	&\mbox{in}\quad\O\tm(0,T),\\
	\wt{u}=0			&\mbox{on}\quad\G\tm(0,T),\\
	\wt{u}(\cdot,0)=u_0		&\mbox{in}\quad\O
\end{cases}
\end{equation}
  and let $u_1\in L^2(\O)$ be the given target function.
  By Theorem 2.1 in \cite{sk}, \eqref{eq:homo} has a unique solution 
$\wt{u}\in C([0,T];L^2(\O))$ and hence $\wt{u}(\cdot,T)\in L^2(\O)$.
  If \eqref{eq:sim} is approximately controllable, then for any $\e>0$ there 
exists $g\in C_0^{\infty}(\G_0\tm(0,T))$ such that the solution $u$ of 
\eqref{eq:sim} satisfies
\[
	\|u(\cdot,T)-\(u_1-\wt{u}(\cdot,T)\)\|_{L^2(\O)}<\e.
\]
  We see that $u+\wt{u}$ solves equation \eqref{eq:ori} and satisfies
\[
	\|(u+\wt{u})(\cdot,T)-u_1\|_{L^2(\O)}<\e.
\]
  Thus approximate controllability for \eqref{eq:sim} immediately implies 
Theorem \ref{thm:ac}.
  In the following, therefore, we will mainly consider \eqref{eq:sim} 
instead of \eqref{eq:ori}.

\section{Representation of the solution}\label{sec:reg}

  In order to obtain the representation of the solution to \eqref{eq:sim}, we 
first prepare the notations.

  Now we are ready to state the following result;

\bigskip
\begin{prop}\label{prop:reg1}
  Let $0<\a<1$ and $g\in C^{\infty}_0(\G_0\tm(0,T))$, then there 
exists a unique solution $u\in C^{\infty}([0,T];H^2(\O))$ to \eqref{eq:sim} 
represented as
\begin{align}\label{eq:gsol}
	u(x,t)=-\sm\(\int_0^t\<g(\cdot,t-\tau),\d_{\nu_L}\v\>
		\tau^{\a-1}E_{\a,\a}(-\l_n\tau^{\a})d\tau\)\v(x).
\end{align}
  The series in \eqref{eq:gsol} is convergent in $C^m([0,T];H^2(\O))$ and 
\begin{align}\label{eq:est02}
	\|\d_t^m u(\cdot,t)\|_{H^2(\O)}
	\le C\(t^{\a(\th-1)+1}\|\d_t^{m+1}g\|_{L^{\infty}(0,T;H^{3/2}(\G))}
		+\|\d_t^m g(\cdot,t)\|_{H^{3/2}(\G)}\)
\end{align}
  for $m=0,1,2,\dots$, where $0<\th<1/4$.
\end{prop}

\bigskip
  In order to prove this proposition, we briefly describe some properties 
concerning convolutions, fractional integrals and the Mittag-Leffler functions.
  First we state the following well known lemma;

\bigskip
\begin{lem}[Young's inequality]\label{lem:young}
  Let $1\le p,q,r\le\infty$ satisfy $1/p+1/q=1+1/r$.
  If $f\in L^p(0,T)$ and $g\in L^q(0,T)$, then the function $f*g$ defined by
\[
	(f*g)(t):=\int_0^t f(t-\tau)g(\tau)d\tau
\]
  belongs to $L^r(0,T)$ and satisfies the estimate
\[
	\|f*g\|_{L^r(0,T)}\le\|f\|_{L^p(0,T)}\|g\|_{L^q(0,T)}.
\]
  In particular, if $r=\infty$, then $f*g$ belongs to $C[0,T]$ (not only 
$L^{\infty}(0,T)$) and
\[
	|(f*g)(t)|\le\|f\|_{L^p(0,t)}\|g\|_{L^q(0,t)},\quad t\in[0,T].
\]
\end{lem}

\bigskip
  For the above lemma, see Appendix A in Stein \cite{st} for example.

  For the convenience of calculation, we introduce the notation of fractional 
integrals.
  For $\a>0$ and $f\in L^1(0,T)$, we define $\a$-th order forward and backward 
integrals of $f$ by
\begin{align*}
	I_{0+}^{\a}f(t):=\frac{1}{\G(\a)}\int_0^t(t-\tau)^{\a-1}f(\tau)d\tau,\\
	I_{T-}^{\a}f(t):=\frac{1}{\G(\a)}\int_t^T(\tau-t)^{\a-1}f(\tau)d\tau.
\end{align*}
  In other words, the forward integral operators of $\a$-th order is the 
convolution with $t^{\a-1}/\G(\a)$ and consequently $I_{0+}^{\a}f$ also 
belongs to $L^1(0,T)$.
  The same argument is also valid for the backward integrals.
  In particular, we have
\begin{equation}\label{eq:cap}
	\d_t^{\a}f(t)=I_{0+}^{1-\a}f'(t)
\end{equation}
  if $0<\a<1$ and $f\in H^1(0,T)$.

  A straightforward calculation yields
\[
	I_{0+}^{\a}\big[t^{\nu}\big]=\frac{\G(\nu+1)}{\G(\nu+\a+1)}t^{\nu+\a}
\]
  for $\nu>-1$ and $\a>0$.
  Therefore, by the termwise integration, we have
\begin{equation}\label{eq:formula}
	I_{0+}^{1-\a}\Big(t^{\a-1}E_{\a,\a}(-\l t^{\a})\Big)
=	E_{\a,1}(-\l t^{\a}),\quad t>0.
\end{equation}
  for $0<\a<1$ and $\l>0$, which is a particular case of (1.100) in \cite{pd}.

  We also have the following formula for fractional integration by parts.

\bigskip
\begin{lem}\label{lem:int}
  Let $\a>0$ and $1<p,q<\infty$ satisfy $1/p+1/q\le1+\a$.
  If $f\in L^p(0,T)$ and $g\in L^q(0,T)$, then
\[
	\int_0^T I_{0+}^{\a}f(t)g(t)dt=\int_0^T f(t)I_{T-}^{\a}g(t)dt.
\]
  In particular, we have
\begin{equation}\label{eq:int}
	(I_{0+}^{\a}f)*g(t)=f*(I_{0+}^{\a}g)(t).
\end{equation}
\end{lem}

\bigskip
  This lemma is derived from Theorem 3.5 in \cite{sm} as its corollary (see 
pp.66-67 in \cite{sm}).

  As for the Mittag-Leffler functions, we have the following two lemmata.

\bigskip
\begin{lem}\label{lem:3.1}
  Let $0<\a<2$ and $\beta\in\R$ be arbitrary and $\mu$ satisfy $\pi\a/2<\mu<\min\{\pi,\pi\a\}$.
  Then there exists a constant $C=C(\a,\beta,\mu)>0$ such that
\begin{equation}\label{eq:3.1}
	|E_{\a,\beta}(z)|\le\frac{C}{1+|z|},\quad\mu\le|\arg(z)|\le\pi.
\end{equation}
\end{lem}

\bigskip
  The proof of Lemma \ref{lem:3.1} can be found on p. 35 in \cite{pd}.

\bigskip
\begin{lem}\label{lem:3.2}
  Let $\l,\a>0$.
  For positive integer $m\in\N$,
\begin{equation}\label{eq:3.2}
	\frac{d^m}{dt^m}E_{\a,1}(-\l t^{\a})
	=-\l t^{\a-m}E_{\a,\a-m+1}(-\l t^{\a}),\quad t>0.
\end{equation}
\end{lem}

\bigskip
\begin{proof}[\bf Proof]
  Since $E_{\a,\beta}(z)$ is an entire function of $z$, the function 
$E_{\a,\beta}(x)$ is real analytic and the series 
$\sum_{k=0}^{\infty}\frac{z^k}{\G(\a k+\beta)}=E_{\a,\beta}(z)$ is termwise 
differentiable in $\R$.
  Since $t^{\a}$ is also real analytic in $t>0$, so is 
$E_{\a,\beta}(-\l t^{\a})$ in $t>0$.
  Therefore equation (\ref{eq:3.2}) can be obtained by termwise 
differentiation.
\end{proof}

\bigskip
  Now we are ready to show Proposition \ref{prop:reg1}.

\bigskip
\begin{proof}[\bf Proof of Proposition \ref{prop:reg1}]
{\bf Step 1.}
  First we prove the unique existence of the solution to \eqref{eq:sim}.
  Since the uniqueness can be shown similarly to Theorem 2.1 in \cite{sk}, it 
is sufficient to show that the solution $u$ of \eqref{eq:sim} is given by 
\eqref{eq:gsol}.

  We split $u$ into $w+\Lm g$ where $w$ solves
\[
\begin{cases}
	\d_t^{\a}w+\L w=-\d_t^{\a}\Lm g	&\mbox{in}\quad\O\tm(0,T),\\
	w=0				&\mbox{on}\quad\G\tm(0,T),\\
	w(\cdot,0)=0			&\mbox{in}\quad\O
\end{cases}
\]
  and
\begin{equation}\label{eq:LgC}
	\Lm g\in C_0^{\infty}((0,T);H^2(\O)).
\end{equation}
  Then $u=w+\Lm g$ satisfies \eqref{eq:sim}.
  By Theorem 2.2 in \cite{sk} (or Proposition 3.1 in \cite{fy}), $w$ is given 
by
\[
	w(x,t)
	=-\sm\(\int_0^t\big((\d_t^{\a}\Lm g)(\cdot,t-\tau),\v\big)
		\tau^{\a-1}E_{\a,\a}(-\l_n\tau^{\a})d\tau\)\v(x).
\]
  Then by equations \eqref{eq:cap}, \eqref{eq:formula} and \eqref{eq:int}, we 
have
\begin{align}\label{eq:u-Lg}
	w(x,t)
&	=-\sm\(\int_0^t\big((I_{0+}^{1-\a}\d_t\Lm g)(\cdot,t-\tau),\v\big)
		\tau^{\a-1}E_{\a,\a}(-\l_n\tau^{\a})d\tau\)\v(x) \nn \\
&	=-\sm\(\int_0^t\big((\d_t\Lm g)(\cdot,t-\tau),\v\big)\cdot
		I_{0+}^{1-\a}\Big(\tau^{\a-1}E_{\a,\a}(-\l_n\tau^{\a})
		\Big)d\tau\)\v(x) \nn \\
&	=-\sm\(\int_0^t\big((\d_t\Lm g)(\cdot,t-\tau),\v\big)
		E_{\a,1}(-\l_n\tau^{\a})d\tau\)\v(x).
\end{align}
  Since $\Lm g(\cdot,0)=0$ by \eqref{eq:LgC}, the integration by parts yields
\begin{align*}
	w(x,t)
&=	\sm\(\int_0^t\frac{\d}{\d\tau}(\Lm g(\cdot,t-\tau),\v)\cdot
			E_{\a,1}(-\l_n\tau^{\a})d\tau\)\v(x) \\
&	=-\Lm g(x,t)-\sm\(\int_0^t(\Lm g(\cdot,t-\tau),\v)\cdot
	\frac{\d}{\d\tau}\Big(E_{\a,1}(-\l_n\tau^{\a})\Big)d\tau\)\v(x).
\end{align*}
  By \eqref{eq:Lgl} and Lemma \ref{lem:3.2}, we have
\begin{align*}
	u(x,t)
&	=w(x,t)+\Lm g(x,t)
	=\sm\(\int_0^t(\Lm g(\cdot,t-\tau),\v)\l_n
		\tau^{\a-1}E_{\a,\a}(-\l_n\tau^{\a})d\tau\)\v(x) \\
&	=-\sm\(\int_0^t\<g(\cdot,t-\tau),\d_{\nu_L}\v\>
		\tau^{\a-1}E_{\a,\a}(-\l_n\tau^{\a})d\tau\)\v(x). \nn
\end{align*}
  Thus the solution $u$ of \eqref{eq:sim} is given by \eqref{eq:gsol}.

\bigskip
\noindent{\bf Step 2.}
  Next we prove that the function $u$ given by (\ref{eq:gsol}) satisfies 
estimate \eqref{eq:est02}.
  Using representation \eqref{eq:u-Lg}, we have
\begin{align*}
	\|Lw(\cdot,t)\|_{L^2(\O)}
&	=\left\|\sm\l_n\(\int_0^t\big((\d_t\Lm g)(\cdot,t-\tau),\v\big)
		E_{\a,1}(-\l_n\tau^{\a})d\tau\)\v\right\|_{L^2(\O)} \\
&	=\left\|\int_0^t\(\sm\l_n\big((\d_t\Lm g)(\cdot,t-\tau),\v\big)
		E_{\a,1}(-\l_n\tau^{\a})\v\)d\tau\right\|_{L^2(\O)} \\
&	\le\int_0^t\left\|\sm\l_n\big((\d_t\Lm g)(\cdot,t-\tau),\v\big)
		E_{\a,1}(-\l_n\tau^{\a})\v\right\|_{L^2(\O)}d\tau \\
&	\le\int_0^t\(\sm\l_n^{2\th}
			\big|\big((\d_t\Lm g)(\cdot,t-\tau),\v\big)\big|^2
		\cdot|\l_n^{1-\th}E_{\a,1}(-\l_n\tau^{\a})|^2\)^{1/2}d\tau.
\end{align*}
  By using Lemma \ref{lem:3.1} again, we have
\[
	|\l_n^{1-\th}E_{\a,1}(-\l_n\tau^{\a})|
\le	\l_n^{1-\th}\cdot\frac{C}{1+\l_n\tau^{\a}}
=	C\cdot\frac{(\l_n\tau^{\a})^{1-\th}}{1+\l_n\tau^{\a}}
		\cdot\tau^{\a(\th-1)}
\le	C\tau^{\a(\th-1)}.
\]
  Therefore,
\begin{align*}
	\|Lw(\cdot,t)\|_{L^2(\O)}
&	\le C\int_0^t\(\sm\l_n^{2\th}\big|\big((\d_t\Lm g)(\cdot,t-\tau),\v\big)
		\big|^2\)^{1/2}\tau^{\a(\th-1)}d\tau \\
&	=C\int_0^t\|(\d_t\Lm g)(\cdot,t-\tau)\|_{\D(L^{\th})}
		\tau^{\a(\th-1)}d\tau \\
&	\le C\int_0^t\|(\d_t g)(\cdot,t-\tau)\|_{H^{3/2}(\G)}
		\tau^{\a(\th-1)}d\tau
\le	C\|\d_t g\|_{L^{\infty}(0,T;H^{3/2}(\G))}
		\int_0^t\tau^{\a(\th-1)}d\tau \\
&	\le Ct^{\a(\th-1)+1}\|\d_t g\|_{L^{\infty}(0,T;H^{3/2}(\G))}.
\end{align*}
  Since $u=w+\Lm g$, we have
\begin{align*}
	\|u(\cdot,t)\|_{H^2(\O)}
&	\le \|w(\cdot,t)\|_{H^2(\O)}+\|\Lm g(\cdot,t)\|_{H^2(\O)}
	\le C\|Lw(\cdot,t)\|_{L^2(\O)}+\|g(\cdot,t)\|_{H^{3/2}(\O)} \\
&	\le C\(t^{\a(\th-1)+1}\|\d_t g\|_{L^{\infty}(0,T;H^{3/2}(\G))}
		+\|g(\cdot,t)\|_{H^{3/2}(\G)}\).
\end{align*}
  Similarly we can also show
\[
	\|\d_t^m u(\cdot,t)\|_{H^2(\O)}
\le	C\(t^{\a(\th-1)+1}\|\d_t^{m+1}g\|_{L^{\infty}(0,T;H^{3/2}(\G))}
		+\|\d_t^m g(\cdot,t)\|_{H^{3/2}(\G)}\)
\]
 for any $m\in\N$.

\bigskip
\noindent{\bf Step 3.}
  We prove that the series in \eqref{eq:gsol} converges in 
$C^m([0,T];H^2(\O))$ for $m=0,1,2,\dots$.
  Since $\Lm g$ clearly belongs to $C^{\infty}([0,T];H^2(\O))$, it is 
sufficient to show the convergence of \eqref{eq:u-Lg}.	
  By the similar calculation to Step 2, we have
\begin{align*}
	&\left\|-\smMN\(\int_0^t\big((\d_t\Lm g)(\cdot,t-\tau),\v\big)
		E_{\a,1}(-\l_n\tau^{\a})d\tau\)\v\right\|_{H^2(\O)} \\
&	\le C\left\|\smMN\(\int_0^t\big((\d_t\Lm g)(\cdot,t-\tau),\v\big)
		E_{\a,1}(-\l_n\tau^{\a})d\tau\)\v\right\|_{\D(L)} \\
&	=C\left\|\smMN\l_n\(\int_0^t\big((\d_t\Lm g)(\cdot,t-\tau),\v\big)
		E_{\a,1}(-\l_n\tau^{\a})d\tau\)\v\right\|_{L^2(\O)} \\
&	\le C\int_0^t\(\smMN\l_n^{2\th}\big|\big((\d_t\Lm g)(\cdot,t-\tau),
		\v\big)\big|^2\)^{1/2}\tau^{\a(\th-1)}d\tau \\
&	\le C\int_0^t\tau^{\a(\th-1)}d\tau\cdot
		\sup_{0\le t\le T}\(\smMN\l_n^{2\th}\big|
		\big((\d_t\Lm g)(\cdot,t),\v\big)\big|^2\)^{1/2} \\
&	\le Ct^{\a(\th-1)+1}\sup_{0\le t\le T}\(\smMN\l_n^{2\th}\big|
		\big((\d_t\Lm g)(\cdot,t),\v\big)\big|^2\)^{1/2}.
\end{align*}
  Since $\d_t\Lm g\in C([0,T];\D(L^{\th}))$, we have
\begin{align*}
	&\sup_{0\le t\le T}
	\left\|-\smMN\(\int_0^t\big((\d_t\Lm g)(\cdot,t-\tau),\v\big)
		E_{\a,1}(-\l_n\tau^{\a})d\tau\)\v\right\|_{H^2(\O)} \\
&	\le CT^{\a(\th-1)+1}\sup_{0\le t\le T}\(\smMN\l_n^{2\th}\big|
		\big((\d_t\Lm g)(\cdot,t),\v\big)\big|^2\)^{1/2}
\to	0\quad\as\ M,N\to\infty.
\end{align*}
  Thus the series in \eqref{eq:u-Lg} is convergent in $H^2(\O)$ uniformly in 
$t\in[0,T]$.
  In the same way, we can also show the uniform convergence of
\[
	\d_t^{m}w(\cdot,t)
=	-\sm\(\int_0^t\big((\d_t^{m+1}\Lm g)(\cdot,t-\tau),\v\big)
	E_{\a,1}(-\l_n\tau^{\a})d\tau\)\v
\]
  for any $m\in\N$.
\end{proof}

\section{Weak solution}\label{sec:weak}

  In this section,

  As we have seen in Proposition \ref{prop:reg1}, the function $u$ defined by 
\eqref{eq:gsol} is the solution of \eqref{eq:ori} with $u_0=0$ when $g$ is 
restricted in $C^{\infty}_0(\G_0\tm(0,T))$.
  However, the domain of the map $g\mapsto u$ can be extended keeping $u$ 
belonging to $C([0,T];L^2(\O))$;

\bigskip
\begin{prop}\label{prop:reg2}
  Let $0<\a<1$ and $g\in L^p(0,T;L^2(\G))$ with $p>4/\a$.
  Then the function $u$ given by \eqref{eq:gsol} belongs to $C([0,T];L^2(\O))$ 
and satisfies
\begin{equation}\label{eq:est}
	\|u(\cdot,t)\|_{L^2(\O)}\le C_2t^{\a\th-1/p}\|g\|_{L^p(0,T;L^2(\G))},
\end{equation}
  where $1/(p\a)<\th<1/4$.
  Moreover for any $0<\de<\th-1/(p\a)$, we have
\begin{equation}\label{eq:Ldu}
	\|u(\cdot,t)\|_{\D(L^{\de})}
	\le Ct^{\a(\th-\de)-1/p}\|g\|_{L^{\infty}(0,T;H^{3/2}(\G))}.
\end{equation}
\end{prop}

\bigskip
\begin{rem}
  If $\a=1$, then the similar result holds for $p>4$ (see \cite{wa}).
\end{rem}

\bigskip
\begin{proof}[\bf Proof of Proposition \ref{prop:reg2}]
\noindent{\bf Step 1.}
  By a simple calculation, we have
\begin{align*}
	\|u(\cdot,t)\|_{L^2(\O)}
&	=\left\|\sm\(\int_0^t\<g(\cdot,t-\tau),\d_{\nu_L}\v\>\tau^{\a-1}
		E_{\a,\a}(-\l_n\tau^{\a})d\tau\)\v\right\|_{L^2(\O)} \\
&	=\left\|\int_0^t\(\sm\<g(\cdot,t-\tau),\d_{\nu_L}\v\>\tau^{\a-1}
		E_{\a,\a}(-\l_n\tau^{\a})\v\)d\tau\right\|_{L^2(\O)} \\
&	\le\int_0^t\left\|\sm\<g(\cdot,t-\tau),\d_{\nu_L}\v\>\tau^{\a-1}
		E_{\a,\a}(-\l_n\tau^{\a})\v\right\|_{L^2(\O)}d\tau \\
&	=\int_0^t\(\sm\left|\<g(\cdot,t-\tau),\d_{\nu_L}\v\>\tau^{\a-1}
		E_{\a,\a}(-\l_n\tau^{\a})\right|^2\)^{1/2}d\tau \\
&	=\int_0^t\(\sm\l_n^{2\th-2}|\<g(\cdot,t-\tau),\d_{\nu_L}\v\>|^2\cdot
		\big|\l_n^{1-\th}\tau^{\a-1}E_{\a,\a}(-\l_n\tau^{\a})\big|^2
		\)^{1/2}d\tau.
\end{align*}
  Similarly to Proposition \ref{prop:reg1}, we use Lemma \ref{lem:3.1} to 
obtain
\[
	|\l_n^{1-\th}\tau^{\a-1}E_{\a,\a}(-\l_n\tau^{\a})|
\le	\l_n^{1-\th}\tau^{\a-1}\cdot\frac{C}{1+\l_n\tau^{\a}}
=	C\cdot\frac{(\l_n\tau^{\a})^{1-\th}}{1+\l_n\tau^{\a}}
		\cdot\tau^{\a\th-1}
\le	C\tau^{\a\th-1}.
\]
  Therefore
\begin{align}\label{eq:g-est}
	\|u(\cdot,t)\|_{L^2(\O)}
&	\le C\int_0^t\(\sm\l_n^{2\th-2}|
		\<g(\cdot,t-\tau),\d_{\nu_L}\v\>|^2\)^{1/2}
		\tau^{\a\th-1}d\tau \nn \\
&	\le C\int_0^t\|g(\cdot,t-\tau)\|_{L^2(\G)}\tau^{\a\th-1}d\tau.
\end{align}
  Let $q\in\R$ satisfy $1/p+1/q=1$, then \eqref{eq:g-est} and Lemma 
\eqref{lem:young} yields
\[
	\|u(\cdot,t)\|_{L^2(\O)}
\le	C\|g\|_{L^p(0,t;L^2(\G))}\(\int_0^t\tau^{q(\a\th-1)}d\tau\)^{1/q}
\le	Ct^{\a\th-1/p}\|g\|_{L^p(0,T;L^2(\G))}.
\]
  Thus we have proved estimate (\ref{eq:est}).
  Moreover, by the similar calculation, we have
\begin{align*}
	&\left\|-\smMN\(\int_0^t\<g(\cdot,t-\tau),\d_{\nu_L}\v\>\tau^{\a-1}
		E_{\a,\a}(-\l_n\tau^{\a})d\tau\)\v\right\|_{L^2(\O)} \\
&	\le C\int_0^t\(\smMN\l_n^{2\th}\big|\big(\Lm g(\cdot,t-\tau),\v\big)
		\big|^2\)^{1/2}\tau^{\a\th-1}d\tau \\
&	\le C\left[\int_0^t\(\smMN\l_n^{2\th}\big|\big(\Lm g(\cdot,\tau),
		\v\big)\big|^2\)^{p/2}d\tau\right]^{1/p}
		\(\int_0^t\tau^{(\a\th-1)q}d\tau\)^{1/q} \\
&	\le Ct^{\a\th-1/p}\left[\int_0^t\(\smMN\l_n^{2\th}\big|\big(
		\Lm g(\cdot,\tau),\v\big)\big|^2\)^{p/2}d\tau\right]^{1/p}.
\end{align*}
  Therefore
\begin{align*}
	&\sup_{0\le t\le T}
	\left\|-\smMN\(\int_0^t\big(\Lm g(\cdot,t-\tau),\v\big)
		E_{\a,1}(-\l_n\tau^{\a})d\tau\)\v\right\|_{H^2(\O)} \\
&	\le CT^{\a\th-1/p}\left[\int_0^T\(\smMN\l_n^{2\th}\big|\big(
		\Lm g(\cdot,\tau),\v\big)\big|^2\)^{p/2}d\tau\right]^{1/p}
\to	0\quad\as\ M,N\to\infty.
\end{align*}
  Thus the series in \eqref{eq:gsol} is convergent in $L^2(\O)$ uniformly in 
$t\in[0,T]$.
  Therefore $u$ belongs to $C([0,T];L^2(\O))$.

\noindent{\bf Step 2.}
  Next we prove \eqref{eq:Ldu}.
  By a simple calculation, we have
\begin{align*}
	\|u(\cdot,t)\|_{\D(L^{\de})}
&	=\left\|\sm\l_n^{\de}\(\int_0^t\<g(\cdot,t-\tau),\d_{\nu_L}\v\>
		\tau^{\a-1}E_{\a,\a}(-\l_n\tau^{\a})d\tau\)\v
		\right\|_{L^2(\O)} \\
&	=\left\|\int_0^t\(\sm\l_n^{\de}\<g(\cdot,t-\tau),\d_{\nu_L}\v\>
		\tau^{\a-1}E_{\a,\a}(-\l_n\tau^{\a})\v\)d\tau
		\right\|_{L^2(\O)} \\
&	\le\int_0^t\left\|\sm\l_n^{\de}\<g(\cdot,t-\tau),\d_{\nu_L}\v\>
		\tau^{\a-1}E_{\a,\a}(-\l_n\tau^{\a})\v\right\|_{L^2(\O)}d\tau\\
&	=\int_0^t\(\sm\left|\l_n^{\de}\<g(\cdot,t-\tau),\d_{\nu_L}\v\>
		\tau^{\a-1}E_{\a,\a}(-\l_n\tau^{\a})\right|^2\)^{1/2}d\tau \\
&	=\int_0^t\(\sm\l_n^{2\th}|(\Lm g(\cdot,t-\tau),\v)|^2\cdot
		\left|\l_n^{1+\de-\th}\tau^{\a-1}E_{\a,\a}(-\l_n\tau^{\a})
			\right|^2\)^{1/2}d\tau.
\end{align*}
  By Lemma \ref{lem:3.1}, we have
\[
	\left|\l_n^{1+\de-\th}\tau^{\a-1}E_{\a,\a}(-\l_n\tau^{\a})\right|
\le	\l_n^{1+\de-\th}\tau^{\a-1}\cdot\frac{C}{1+\l_n\tau^{\a}}
=	C\cdot\frac{(\l_n\tau^{\a})^{1+\de-\th}}{1+\l_n\tau^{\a}}
		\cdot\tau^{\a(\th-\de)-1}
\le	C\tau^{\a(\th-\de)-1}.
\]
  Therefore combining this with \eqref{eq:Lg2}, we have
\begin{align*}
	\|u(\cdot,t)\|_{\D(L^{\de})}
&	\le C\int_0^t\(\sm\l_n^{2\th}|(\Lm g(\cdot,t-\tau),\v)|^2\)^{1/2}
		\tau^{\a(\th-\de)-1}d\tau \\
&	=C\int_0^t\|\Lm g(\cdot,t-\tau)\|_{\D(L^{\th})}
		\tau^{\a(\th-\de)-1}d\tau
	\le C\int_0^t\|g(\cdot,t-\tau)\|_{L^2(\G)}
		\tau^{\a(\th-\de)-1}d\tau.
\end{align*}
  By taking $q\in\R$ as before, we apply Lemma \ref{lem:young} again and have
\begin{align*}
	\|u(\cdot,t)\|_{\D(L^{\de})}
\le	C\|g\|_{L^p(0,t;L^2(\G))}
		\(\int_0^t\tau^{q(\a(\th-\de)-1)}d\tau\)^{1/q}
\le	Ct^{\a(\th-\de)-1/p}\|g\|_{L^p(0,T;L^2(\G))}.
\end{align*}
  Thus the proof of Proposition \ref{prop:reg2} is completed.
\end{proof}

\section{Dual system}\label{sec:dual}

  We prove the following propositions;

\bigskip
\begin{prop}\label{pr:ana}
  Let $0<\a<1$ and $v_0\in L^2(\O)$.
  Then there exists a unique solution $v\in C([0,T);H^2(\O)\cap H^1_0(\O))$ to 
$(\ref{eq:dual})$ which is represented as
\begin{equation}\label{eq:dsol2}
	v(x,t)=\sm(v_0,\v)(T-t)^{\a-1}E_{\a,\a}(-\l_n(T-t)^{\a})\v(x).
\end{equation}
  and has the following estimate for any $0<\de\le1$;
\begin{equation}\label{eq:Ldv}
	\|v(\cdot,t)\|_{\D(L^{1-\de})}\le C(T-t)^{\a\de-1}\|v_0\|_{L^2(\O)}.
\end{equation}
  Moreover the mapping $[0,T)\ni t\mapsto\d_{\nu_L}v(\cdot,t)\in L^2(\G)$ is 
analytically extended to $S_T:=\{z\in\C;\re z<T\}$.
\end{prop}

\bigskip
\begin{prop}\label{pr:unicon}
  Let $\G_0$ be open in $\G$ and $v\in C([0,T);H^2(\O)\cap H^1_0(\O))$ be the 
solution of \eqref{eq:dual} corresponding to $v_0\in L^2(\O)$.
  If $\d_{\nu_L}v=0$ on $\G_0\tm(0,T)$, then $v=0$ in $\O\tm(0,T)$.
\end{prop}

\bigskip
\begin{proof}[\bf Proof of Proposition \ref{pr:ana}]
  By Proposition 4.1 in \cite{fy}, it is already known that \eqref{eq:dual} 
has a unique solution and that it is given by \eqref{eq:dsol2}.

  We first show estimate \eqref{eq:Ldv}.
  By \eqref{eq:dsol2}, we have
\begin{align*}
	\|v(\cdot,t)\|_{\D(L^{1-\de})}^2
&	=\left\|\sm\l_n^{1-\de}(v_0,\v)(T-t)^{\a-1}E_{\a,\a}(-\l_n(T-t)^{\a})
		\v\right\|_{L^2(\O)}^2 \\
&	=\sm|\l_n^{1-\de}(v_0,\v)(T-t)^{\a-1}E_{\a,\a}(-\l_n(T-t)^{\a})|^2.
\end{align*}
  We use Lemma \ref{lem:3.1} to obtain
\begin{align*}
	|\l_n^{1-\de}(T-t)^{\a-1}E_{\a,\a}(-\l_n(T-t)^{\a})|
&	\le\l_n^{1-\de}(T-t)^{\a-1}\cdot\frac{C}{1+\l_n(T-t)^{\a}} \\
&	=C\cdot\frac{\(\l_n(T-t)^{\a}\)^{1-\de}}{1+\l_n(T-t)^{\a}}
		\cdot(T-t)^{\a\de-1} \\
&	\le C(T-t)^{\a\de-1}.
\end{align*}
  Therefore,
\[
	\|v(\cdot,t)\|_{\D(L^{1-\de})}
	\le C(T-t)^{\a\de-1}\(\sm|(v_0,\v)|^2\)^{1/2}
	=C(T-t)^{\a\de-1}\|v_0\|_{L^2(\O)}.
\]

  Next we show the analyticity of $\d_{\nu_L}v(\cdot,t)$ in $t\in S_T$.
  Since $\d_{\nu_L}:H^2(\O)\to L^2(\G)$ is bounded, we have
\begin{equation}\label{eq:neu}
	\frac{\d v}{\d\nu_L}(\cdot,t)
	=\sm(v_0,\v)(T-t)^{\a-1}E_{\a,\a}(-\l_n(T-t)^{\a})\frac{\d\v}{\d\nu_L}
\end{equation}
  and the right-hand side of the above is convergent in $L^2(\G)$ for any 
$t\in(0,T)$.

  We note that $E_{\a,\a}(-\l_n z)$ is an entire function (see Section 1.8 in 
\cite{ks} for example) and consequently each 
$(T-z)^{\a-1}E_{\a,\a}(-\l_n(T-z)^\a)$ is analytic in $z\in S_T$.
  Therefore 
$S_T\ni z\mapsto
\smN(v_0,\v)(T-z)^{\a-1}E_{\a,\a}(-\l_n(T-z)^\a)\d_{\nu_L}\v\in L^2(\G)$ 
is also analytic.
  If we fix $\de'>0$ arbitrarily, then for $z\in\C$ with $\re z\le T-\de'$, we 
have
\begin{align*}
	&\left\|\smMN(v_0,\v)(T-z)^{\a-1}E_{\a,\a}(-\l_n(T-z)^\a)\d_{\nu_L}\v
		\right\|^2_{L^2(\G)} \\
&	\le C\left\|\smMN(v_0,\v)(T-z)^{\a-1}E_{\a,\a}(-\l_n(T-z)^\a)\v
		\right\|^2_{H^2(\O)} \\
&	\le C\left\|\smMN(v_0,\v)(T-z)^{\a-1}E_{\a,\a}(-\l_n(T-z)^\a)\v
		\right\|^2_{\D(L)} \\
&	=C\smMN|\l_n(v_0,\v)(T-z)^{\a-1}E_{\a,\a}(-\l_n(T-z)^\a)|^2 \\
&	\le C\smMN|(v_0,\v)|^2|T-z|^{-2}
		\(\frac{\l_n|T-z|^\a}{1+\l_n|T-z|^\a}\)^2 \\
&	\le C\de'^{-2}\smMN|(v_0,\v)|^2
	\to 0 \quad\as\ M,N\to\infty.
\end{align*}
  That is, \eqref{eq:neu} is uniformly convergent in 
$\{z\in\C;\re z\le T-\de'\}$.
  Hence $\d_{\nu_L}v(\cdot,t)$ is also analytic in $t\in S_T$.
\end{proof}

\bigskip
\begin{proof}[\bf Proof of Proposition \ref{pr:unicon}]
  Since $\d_{\nu_L}v=0$ in $\G_0\tm(0,T)$ and 
$\d_{\nu_L}v:[0,T)\to L^2(\G)$ can be analytically extended to $S_T$, we have
\begin{equation}\label{eq:4.6}
	\frac{\d v}{\d\nu_L}(x,t)
=	\sm(v_0,\v)(T-t)^{\a-1}E_{\a,\a}(-\l_n(T-t)^{\a})
	\frac{\d\v}{\d\nu_L}(x)
=	0,	\quad x\in\G_0,\ t\in(-\infty,T).
\end{equation}
  Let $\{\mu_k\}_{k\in\N}$ be all spectra of $L$ without multiplicities and we 
denote by $\{\vk\}_{1\le j\le m_k}$ an orthonormal basis of $\kr(\mu_k-L)$.
  By using these notations, we can rewrite (\ref{eq:4.6}) by
\begin{equation}\label{eq:4.7}
	\smk\(\smj(v_0,\vk)\frac{\d \vk}{\d\nu_L}(x)\)(T-t)^{\a-1}
		E_{\a,\a}(-\mu_k(T-t)^\a)
=	0,\quad x\in\G_0,\ t\in(-\infty,T).
\end{equation}
  We regard $\d_{\nu_L}$ as a bounded operator from $H^{2\e+3/2}(\O)$ to 
$H^{2\e}(\G)$ with $0<\e<1/4$.
  Then for any $z\in\C$ with $\re z=\xi>0$ and $N\in\N$, we have
\begin{align*}
&	\left\|\sum_{k=1}^N\(\smj(v_0,\vk)\frac{\d \vk}{\d\nu_L}\)e^{z(t-T)}
		(T-t)^{\a-1}E_{\a,\a}(-\mu_k(T-t)^\a)\right\|^2_{L^2(\G)} \\
&	\le\left\|\sum_{k=1}^N\(\smj(v_0,\vk)\frac{\d \vk}{\d\nu_L}\)e^{z(t-T)}
		(T-t)^{\a-1}E_{\a,\a}(-\mu_k(T-t)^\a)\right\|^2_{H^{2\e}(\G)}\\
&	\le C\left\|\sum_{k=1}^N\(\smj(v_0,\vk)\vk\)e^{z(t-T)}(T-t)^{\a-1}
		E_{\a,\a}(-\mu_k(T-t)^\a)\right\|^2_{H^{2\e+3/2}(\O)} \\
&	\le C\left\|\sum_{k=1}^N\(\smj(v_0,\vk)\vk\)e^{z(t-T)}(T-t)^{\a-1}
		E_{\a,\a}(-\mu_k(T-t)^\a)\right\|^2_{\D(L^{\e+3/4})} \\
&	=C\sum_{k=1}^N\(\smj|(v_0,\vk)|^2\)e^{2\xi(t-T)}\left|\mu_k^{\e+3/4}
		(T-t)^{\a-1}E_{\a,\a}(-\mu_k(T-t)^\a)\right|^2.
\end{align*}
  By Lemma \ref{lem:3.1}, we have
\begin{align*}
	\left|\mu_k^{\e+3/4}(T-t)^{\a-1}E_{\a,\a}(-\mu_k(T-t)^\a)\right|
&\le	\mu_k^{\e+3/4}(T-t)^{\a-1}\cdot\frac{C}{1+\mu_k(T-t)^\a} \\
&\le	C(T-t)^{\a(1/4-\e)-1}
		\cdot\frac{\(\mu_k(T-t)^\a\)^{\e+3/4}}{1+\mu_k(T-t)^\a} \\
&\le	C(T-t)^{\a(1/4-\e)-1}.
\end{align*}
  Therefore
\begin{align*}
	&\left\|\sum_{k=1}^N\(\smj(v_0,\vk)\d_{\nu_L}\vk\)e^{z(t-T)}
		(T-t)^{\a-1}E_{\a,\a}(-\mu_k(T-t)^\a)\right\|_{L^2(\G)} \\
&	\le Ce^{\xi(t-T)}(T-t)^{\b-1}\|v_0\|_{L^2(\O)}.
\end{align*}
  where $\b:=(1/4-\e)\a>0$.
  The right-hand side of the above is integrable on $(-\infty,T)$;
\[
	\int_{-\infty}^T e^{\xi(t-T)}(T-t)^{\b-1}dt
	=\frac{\G(\b)}{\xi^{\b}}.
\]
  Hence the Lebesgue theorem yields that
\begin{flalign}\label{eq:4.75}
	&\int_{-\infty}^T e^{z(t-T)}\(\smk\smj(v_0,\vk)
		\frac{\d \vk}{\d\nu_L}(x)(T-t)^{\a-1}
		E_{\a,\a}(-\mu_k(T-t)^\a)\)dt \nn \\
&	=\smk\smj(v_0,\vk)\frac{\d \vk}{\d\nu_L}(x)
		\(\int_{-\infty}^T e^{z(t-T)}
		(T-t)^{\a-1}E_{\a,\a}(-\mu_k(T-t)^\a)dt\) \nn \\
&	=\smk\smj\frac{(v_0,\vk)}{z^\a+\mu_k}\frac{\d \vk}{\d\nu_L}(x),
		\quad\al x\in\G,\ \re z>0,
\end{flalign}
  where we have used the Laplace transform formula;
\[
	\int_0^{\infty}e^{-zt}t^{\a-1}E_{\a,\a}(-\mu_k t^{\a})dt
=	\frac{1}{z^{\a}+\mu_k}, \quad\re z>0
\]
  (see (1.80) in p.21 of \cite{pd}).
  By (\ref{eq:4.7}) and (\ref{eq:4.75}), we have
\[
	\smk\smj\frac{(v_0,\vk)}{z^\a+\mu_k}\frac{\d \vk}{\d\nu_L}(x)=0,
	\quad \al x\in\G_0,\ \re z>0,
\]
  that is,
\[
	\smk\smj\frac{(v_0,\vk)}{\eta+\mu_k}\frac{\d \vk}{\d\nu_L}(x)=0,
	\quad \al x\in\G_0,\ \re\eta>0.
\]
  By using analytic continuation in $\eta$, we have
\begin{equation}\label{eq:4.10}
	\smk\smj\frac{(v_0,\vk)}{\eta+\mu_k}\frac{\d \vk}{\d\nu_L}(x)=0,
	\quad \al x\in\G_0,\ \eta\in\C\setminus\{-\mu_k\}_{k\in\N}.
\end{equation}
  Then we can take a suitable disk which includes $-\mu_\ell$ and does not 
include $\{-\mu_k\}_{k\ne\ell}$.
  By integrating (\ref{eq:4.10}) in the disk, we have
\[
	\smjl(v_0,\vlj)\frac{\d\vlj}{\d\nu_L}(x)=0,	\quad \al x\in\G_0.
\]
  By setting 
$\wt{v}_{\ell}:=\smjl(v_0,\vlj)\vlj$, 
we have
\[
	(L-\mu_{\ell})\wt{v}_\ell=0	\quad\mbox{in $\O$}	\quad\an\quad
	\frac{\d\wt{v}_\ell}{\d\nu_L}=0	\quad\mbox{on $\G_0$}.
\]
  Therefore the unique continuation result for eigenvalue problem of elliptic 
operator (see Corollary 2.2 in \cite{sw} or Chapter 3 in \cite{is} for 
example) implies
\[
	\wt{v}_\ell(x)=\smjl(v_0,\vlj)\vlj(x)=0,	\quad x\in\O
\]
  for each $\ell\in\N$.
  Since $\{\vlj\}_{1\le j\le m_\ell}$ is linearly independent in 
$\O$, we see that
\[
	(v_0,\vlj)=0,	\quad 1\le j\le m_\ell,\ \ell\in\N.
\]
  This implies $v=0$ in $\O\tm(0,T)$.
\end{proof}

\section{Proof of main results}\label{sec:proof}

  In this section, we complete the proof of our main theorems.

\bigskip
\begin{proof}[\bf Proof of Theorem \ref{thm:ac}]
{\bf Step 1.}
  We first show that for any $g\in C_0^{\infty}(\G_0\times(0,T))$ and 
$v_0\in L^2(\O)$, the following identity holds;
\begin{equation}\label{eq:Fta}
	\int_{\O}u(x,T)v_0(x)dx
	+\int_0^T\int_{\G_0}g(x,t)\frac{\d v}{\d\nu_L}(x,t)\ds_xdt=0,
\end{equation}
  where $u$ and $v$ are the corresponding solutions of \eqref{eq:sim} and 
\eqref{eq:dual} respectively.
  Since the first equation in \eqref{eq:sim} holds in 
$C^{\infty}([0,T];L^2(\O))$ by Proposition \ref{prop:reg1} and 
$v\in L^1(0,T;L^2(\O))$ by \eqref{eq:Ldv} with $\de=1$, we see that
\[
	0=\int_0^T\int_{\O}\big(\cd^{\a}u+\L u\big)v dxdt
	=\int_0^T\int_{\O}(\cd^{\a}u)v dxdt+\int_0^T\int_{\O}(\L u)v dxdt.
\]
  In the above equation, the first term is calculated as follows;
\begin{align*}
	\int_0^T\int_{\O}(\cd^{\a}u)vdxdt
&	=\int_0^T\int_{\O}I_{0+}^{1-\a}\frac{\d u}{\d t}\cdot v dxdt
	=\int_0^T\int_{\O}\frac{\d u}{\d t}\cdot I_{T-}^{1-\a}v dxdt\\
&	=\left.\int_{\O}u\(I_{T-}^{1-\a}v\)dx\right|_{t=0}^{t=T}
		-\int_0^T\int_{\O}u\cdot\frac{\d}{\d t}I_{T-}^{1-\a}v dxdt \\
&	=\int_{\O}u(\cdot,T)v_0 dx+\int_0^T\int_{\O}u(D_t^{\a}v)dxdt.
\end{align*}
  Here we have used Lemma \ref{lem:int}, the integration in $t$ by parts and 
the initial conditions in \eqref{eq:sim} and \eqref{eq:dual}.
  In terms of $u\in C^{\infty}([0,T];H^2(\O))$ and 
$v \in C([0,T);H^2(\O)\cap H^1_0(\O))$ by Propositions \ref{prop:reg1} and 
\ref{pr:ana}, we apply the Green formula to have
\begin{align*}
	\int_0^T\int_{\O}(\L u)v dxdt
&	=\int_0^T\int_{\O}u(\L v)dxdt
		+\int_0^T\int_{\G}\(u\frac{\d v}{\d\nu_L}
			-\frac{\d u}{\d\nu_L}v\)\ds_x dt \\
&	=\int_0^T\int_{\O}u(\L v)dxdt
		+\int_0^T\int_{\G_0}g\frac{\d v}{\d\nu_L}\ds_x dt.
\end{align*}
  In the above calculation, we have used boundary conditions in \eqref{eq:sim} 
and \eqref{eq:dual}.
  We also note that by \eqref{eq:Ldu} and \eqref{eq:Ldv}, the function 
\[
	t\mapsto\int_{\O}u(x,t)\L v(x,t)dx=(u(\cdot,t),Lv(\cdot,t))
\]
  is integrable in $t\in(0,T)$.
  Therefore we have
\begin{align*}\label{eq:var}
0&	=\int_0^T\int_{\O}(\cd^{\a}u)v dxdt+\int_0^T\int_{\O}(\L u)v dxdt\nn \\
&	=\left\{\int_{\O}u(\cdot,T)v_0 dx
		+\int_0^T\int_{\O}u(D_t^{\a}v)dxdt\right\}
	+\left\{\int_0^T\int_{\O}u(\L v)dxdt
		+\int_0^T\int_{\G_0}g\frac{\d v}{\d\nu_L}\ds_x dt\right\} \nn\\
&	=\int_{\O}u(\cdot,T)v_0 dx
		+\int_0^T\int_{\O}u\big(D_t^{\a}v+\L v\big)dxdt
		+\int_0^T\int_{\G_0}g\frac{\d v}{\d\nu_L}\ds_x dt \nn \\
&	=\int_{\O}u(\cdot,T)v_0 dx
		+\int_0^T\int_{\G_0}g\frac{\d v}{\d\nu_L}\ds_xdt.
\end{align*}
  Thus we have proved \eqref{eq:Fta}.

\noindent{\bf Step 2.}
  We note that the assertion of Theorem \ref{thm:ac} is equivalent to
\begin{equation}\label{eq:perp}
	\{u(\cdot,T);\ g\in C_0^{\infty}(\G_0\times(0,T))\}^{\perp}=\{0\},
\end{equation}
  where the orthogonal complement is taken in $L^2(\O)$.
  Suppose that $v_0\in L^2(\O)$ satisfies
\[
	(u(\cdot,T),v_0)=0
\]
  for any $g\in C_0^{\infty}(\G_0\tm(0,T))$.
  Then, by \eqref{eq:Fta}, we have
\[
	\int_0^T\int_{\G_0}g(x,t)\frac{\d v}{\d\nu_L}(x,t)\ds_x dt=0
\]
  for any $g\in C_0^{\infty}(\G_0\tm(0,T))$.
  By the fundamental lemma of the calculus of variations, we have
\[
	\frac{\d v}{\d\nu_L}(x,t)=0, \quad(x,t)\in\G_0\tm(0,T),
\]
  from which Proposition \ref{pr:unicon} implies
\[
	v_0\equiv0.
\]
  Thus we have shown \eqref{eq:perp} and completed the proof of Theorem 
\ref{thm:ac}.
\end{proof}

\bigskip
\begin{proof}[\bf Proof of Theorem \ref{thm:ch}]
{\bf Step 1.}
  First we show that $J_{\e}$ admits a unique minimizer.
  Since $J_{\e}$ is clearly convex and lower semi-continuous, it suffices to 
show its coercivity.

  Let $\{v_0^j\}$ be a sequence in $L^2(\O)$ such that
\[
	\lim_{j\to\infty}\|v_0^j\|_{L^2(\O)}=\infty.
\]
  We set $w_0^j:=v_0^j/\|v_0^j\|_{L^2(\O)}$ and denote by $w^j$ the solution 
of \eqref{eq:dual} with $v_0=w_0^j$.
  That is,
\begin{equation}\label{eq:wj}
	w^j(x,t)=\sm(w_0^j,\v)(T-t)^{\a-1}E_{\a,\a}(-\l_n(T-t)^{\a})\v(x)
\end{equation}
  (see \eqref{eq:dsol2}).
  Then we have
\begin{align*}
	\frac{J_{\e}(v_0^j)}{\|v_0^j\|_{L^2(\O)}}
=	\frac{\|v_0^j\|_{L^2(\O)}}{2}\int_0^T\int_{\G_0}(T-t)^2
		|\d_{\nu_L}w^j|^2\ds_xdt
	+\e\|(I-\pi_E)w_0^j\|_{L^2(\O)}+(w_0^j,u_1).
\end{align*}
  If 
$\li_{j\to\infty}\int_0^T\int_{\G_0}(T-t)^2|\d_{\nu_L}w^j|^2\ds_xdt>0$, 
then we immediately obtain
\[
	\lim_{j\to\infty}\frac{J_{\e}(v_0^j)}{\|v_0^j\|_{L^2(\O)}}=\infty,
\]
  from which the coercivity of $J_{\e}$ follows.
  In the following, therefore, we assume
\begin{equation}\label{eq:lim}
	\li_{j\to\infty}
	\int_0^T\int_{\G_0}(T-t)^2|\d_{\nu_L}w^j|^2\ds_xdt=0.
\end{equation}
  Since $\|w_0^j\|_{L^2(\O)}=1$, there exists a subsequence (denoted by 
$\{w_0^j\}$ again without any confusion) weakly converging to some 
$\ov{w}_0\in L^2(\O)$.
  Let $\ov{w}$ be the solution of \eqref{eq:dual} wth $v_0=\ov{w}_0$, that is, 
\begin{equation}\label{eq:w}
	\ov{w}(x,t)
	=\sm(\ov{w}_0,\v)(T-t)^{\a-1}E_{\a,\a}(-\l_n(T-t)^{\a})\v(x).
\end{equation}
  Then we see that
\begin{equation}\label{eq:wjw}
	\frac{\d w^j}{\d\nu_L}(\cdot,t)\to\frac{\d\ov{w}}{\d\nu_L}(\cdot,t)
\end{equation}
  in $L^2(\G)$ for any $t\in (0,T)$.
  Indeed, by \eqref{eq:wj} and \eqref{eq:w}, we use the similar calculation 
to the proof of Proposition \ref{pr:ana} and have
\begin{align*}
	\left\|\frac{\d w^j}{\d\nu_L}(\cdot,t)
		-\frac{\d\ov{w}}{\d\nu_L}(\cdot,t)\right\|_{L^2(\G)}^2
&\le 	C\left\|w^j(\cdot,t)-\ov{w}(\cdot,t)\right\|_{H^2(\O)}^2
\le	C\left\|w^j(\cdot,t)-\ov{w}(\cdot,t)\right\|_{\D(L)}^2 \\
&=	C\sm\l_n^2|(w_0^j-\ov{w}_0,\v)|^2(T-t)^{2\a-2}
		|E_{\a,\a}(-\l_n(T-t)^{\a})|^2 \\
&\le	C\sm|(w_0^j-\ov{w}_0,\v)|^2(T-t)^{-2}
		\left|\frac{\l_n(T-t)^{\a}}{1+\l_n(T-t)^{\a}}\right|^2 \\
&\le	\frac{C}{(T-t)^2}\sm|(w_0^j-\ov{w}_0,\v)|^2
\to	0
\end{align*}
  for any $t\in(0,T)$, from which \eqref{eq:wjw} follows.
  Here we have used Lebesgue's convergence theorem regarding the summation 
as an integral.
  Now we set
\[
	\psi^j(t):=(T-t)^2\left\|
		\frac{\d w^j}{\d\nu_L}(\cdot,t)\right\|_{L^2(\G)}^2
	\quad\an\quad
	\ov{\psi}(t):=(T-t)^2\left\|
		\frac{\d\ov{w}}{\d\nu_L}(\cdot,t)\right\|_{L^2(\G)}^2.
\]
  Then by \eqref{eq:wjw},
\[
	\lim_{j\to\infty}\psi^j(t)=\ov{\psi}(t),\quad t\in(0,T).
\]
  Moreover, by the representation of \eqref{eq:wj}, we have
\begin{align*}
	\psi^j(t)
&=	(T-t)^2\left\|\frac{\d w^j}{\d\nu_L}(\cdot,t)\right\|_{L^2(\G)}^2 
\le 	C(T-t)^2\left\|w^j(\cdot,t)\right\|_{H^2(\O)}^2
\le	C(T-t)^2\left\|w^j(\cdot,t)\right\|_{\D(L)}^2 \\
&=	C(T-t)^2\sm\l_n^2|(w_0^j,\v)|^2(T-t)^{2\a-2}
		|E_{\a,\a}(-\l_n(T-t)^{\a})|^2 \\
&\le	C\sm|(w_0^j,\v)|^2
		\left|\frac{\l_n(T-t)^{\a}}{1+\l_n(T-t)^{\a}}\right|^2
\le	C\|w_0^j\|_{L^2(\G)}^2=C.
\end{align*}
  Therefore, by Lebesgue's convergence theorem, we have
\[
	\int_0^T\psi^j(t)dt\to\int_0^T\ov{\psi}(t)dt,
\]
  that is, 
\[
	\int_0^T\int_{\G_0}(T-t)^2|\d_{\nu_L}w^j|^2\ds_xdt
	\to\int_0^T\int_{\G_0}(T-t)^2|\d_{\nu_L}\ov{w}|^2\ds_xdt.
\]
  Combining this with \eqref{eq:lim}, we have
\[
	\int_0^T\int_{\G_0}(T-t)^2|\d_{\nu_L}\ov{w}|^2\ds_xdt=0.
\]
  Hence we have
\[
	\d_{\nu_L}\ov{w}=0\quad\mbox{on}\ \G_0\tm(0,T),
\]
  from which we deduce $\ov{w}_0\equiv0$ by Proposition \ref{pr:unicon}.
  That is, $\{w_0^j\}$ weakly converges to 0 in $L^2(\O)$ and consequently we 
have
\[
	\lim_{j\to\infty}(w_0^j,u_1)=0
	\quad\an\quad
	\lim_{j\to\infty}\|(I-\pi_E)w_0^j\|_{L^2(\O)}=1
\]
  since $\pi_E$ is a compact operator.
  Therefore we obtain
\[
	\li_{j\to\infty}\frac{J_{\e}(v_0^j)}{\|v_0^j\|_{L^2(\O)}}\ge\e.
\]
  Thus we have shown the coercivity of $J_{\e}$.

\noindent{\bf Step 2.}
  Let $\ov{v}_0$ be the minimizer of $J_{\e}$, then for any $h>0$ and 
$v_0\in L^2(\O)$, we have
\begin{align*}
0&\le	J_{\e}(\ov{v}_0+hv_0)-J_{\e}(\ov{v}_0) \\
&\le	h\int_0^T\int_{\G_0}(T-t)^2(\d_{\nu_L}\ov{v})(\d_{\nu_L}v)\ds_xdt
	+\frac{h^2}{2}\int_0^T\int_{\G_0}(T-t)^2|\d_{\nu_L}v|^2\ds_xdt \\
&\qquad	+\e h\|(I-\pi_E)v_0\|_{L^2(\O)}+h(v_0,u_1).
\end{align*}
  Dividing the above inequality by $h$ and letting $h\to0$, we have
\begin{align*}
0&\le	\int_0^T\int_{\G_0}(T-t)^2(\d_{\nu_L}\ov{v})(\d_{\nu_L}v)\ds_xdt
		+\e\|(I-\pi_E)v_0\|_{L^2(\O)}-(v_0,u_1) \\
&=	\int_0^T\int_{\G_0}g(\d_{\nu_L}v)\ds_xdt+\e\|(I-\pi_E)v_0\|_{L^2(\O)}
		+(v_0,u_1).
\end{align*}
  By the density argument, we can verify \eqref{eq:Fta} for 
$g\in L^p(0,T;L^2(\O))$ with $p>4/\a$.
  Then we obtain
\[
	0\le-(u(\cdot,T),v_0)+\e\|(I-\pi_E)v_0\|_{L^2(\O)}+(u_1,v_0),
\]
  that is,
\[
	(u(\cdot,T)-u_1,v_0)\le\e\|(I-\pi_E)v_0\|_{L^2(\O)}.
\]
  By taking $h<0$ and repeating the same argument, we also have 
$(u_1-u(\cdot,T),v_0)\le\e\|(I-\pi_E)v_0\|_{L^2(\O)}$.
  Therefore
\begin{equation}\label{eq:v0}
	|(u(\cdot,T)-u_1,v_0)|\le\e\|(I-\pi_E)v_0\|_{L^2(\O)}.
\end{equation}
  Since $v_0\in L^2(\O)$ was arbitrary, we take $v_0\in E^{\perp}$ and 
obtain
\[
	|(u(\cdot,T)-u_1,v_0)|\le\e\|v_0\|_{L^2(\O)},
\]
  that is,
\[
	\|u(\cdot,T)-u_1\|_{L^2(\O)}\le\e.
\]
  Moreover, by taking $v_0\in E$ in \eqref{eq:v0}, we have
\[
	|(u(\cdot,T)-u_1,v_0)|=0.
\]
  Since $v_0\in E$ can be taken arbitrarily, we have
\[
	u(\cdot,T)-u_1\in E^{\perp},
\]
  that is,
\[
	\pi_E(u(\cdot,T))=\pi_E(u_1).
\]
\end{proof}

\appendix
\section{Regularity of the elliptic problem}

  In this section, we consider the following elliptic boundary value problem;
\begin{equation}\label{eq:ell}
\begin{cases}
	\L u=0	&\mbox{in}\quad\O, \\
	u=g	&\mbox{on}\quad\G,
\end{cases}
\end{equation}
  where $g$ is given on $\G$.
  For $g\in H^{3/2}(\G)$, by using the trace theorem and lifting and applying 
the well known results for the elliptic boudndary value problems with 
homogeneous data (see Theorems 8.1 and 9.8 in Agmon \cite{am} for example), we 
see that \eqref{eq:ell} has a unique solution $u\in H^2(\O)$ satisfying
\begin{equation}\label{eq:ug}
	\|u\|_{H^2(\O)}\le C\|g\|_{H^{3/2}(\G)}.
\end{equation}
  In the following, we will discuss \eqref{eq:ell} for non-smooth $g$ by the 
{\it transposition method}.

  We first consider the {\it dual system} for \eqref{eq:ell};
\begin{equation}\label{eq:ell'}
\begin{cases}
	\L v=f	&\mbox{in}\quad\O, \\
	v=0	&\mbox{on}\quad\G,
\end{cases}
\end{equation}
  where $f$ is given in $\O$.
  It is well known that for any $f\in L^2(\O)$, \eqref{eq:ell'} posesses a 
unique solution $v\in H^2(\O)$ satisfying
\begin{equation}\label{eq:vf}
	\|v\|_{H^2(\O)}\le C\|f\|_{L^2(\O)}.
\end{equation}
  Henceforth we will denote this solution by $v_f$.
  Now we can define the solution of \eqref{eq:ell} in a weaker sense.

\bigskip
\begin{df}
  A function $u$ is a {\it weak solution} of \eqref{eq:ell} if 
\begin{equation}\label{eq:weak}
	\int_{\O}u(x)f(x)dx+\int_{\G}g(x)\frac{\d v_f}{\d\nu_L}(x)\ds_x=0
\end{equation}
  holds for any $f\in L^2(\O)$.
\end{df}

\bigskip
  According to this definition, the solution $u\in H^2(\O)$ obtained before is 
also a weak solution.
  Indeed, by the Green's formula, we have
\begin{align*}
0	&=\int_{\O}\L u(x)v_f(x)dx
	=\int_{\O}u(x)\L v_f(x)dx
	+\int_{\G}\(u(x)\frac{\d v_f}{\d\nu_L}(x)
		-\frac{\d u}{\d\nu_L}(x)v(x)\)\ds_x \\
	&=\int_{\O}u(x)f(x)dx+\int_{\G}g(x)\frac{\d v_f}{\d\nu_L}(x)\ds_x.
\end{align*}
  Thus condition \eqref{eq:weak} is satisfied.
  We also see that \eqref{eq:ell} has a weak solution if $g$ is a distribution;

\bigskip
\begin{prop}\label{prop:weak}
  For any $g\in H^{-1/2}(\G)$, there exists a unique weak solution 
$u\in L^2(\O)$ satisfying
\begin{equation}\label{eq:ug'}
	\|u\|_{L^2(\O)}\le C\|g\|_{H^{-1/2}(\G)}.
\end{equation}
\end{prop}

\bigskip
\begin{proof}
  As we have seen, for $f\in L^2(\O)$, the solution $v_f$ of \eqref{eq:ell'} 
belongs to $H^2(\O)$.
  Therefore, by the trace theorem, $\d_{\nu_L}v_f\in H^{1/2}(\G)$ and 
\[
	\left\|\frac{\d v_f}{\d\nu_L}\right\|_{H^{1/2}(\G)}
\le	C\|v_f\|_{H^2(\O)}.
\]
  Combining this with \eqref{eq:vf}, we obtain
\begin{equation}\label{eq:dvf}
	\left\|\frac{\d v_f}{\d\nu_L}\right\|_{H^{1/2}(\G)}
\le	C\|f\|_{L^2(\O)}.
\end{equation}
  Thus the mapping
\[
	L^2(\O)\ni f\mapsto\frac{\d v_f}{\d\nu_L}\in H^{1/2}(\G)
\]
  is bounded, and so is 
\[
	L^2(\O)\ni f\mapsto-\int_{\G}g(x)\frac{\d v_f}{\d\nu_L}(x)\ds_x\in\C.
\]
  Hence the Riesz's representation theorem yields that there exists a unique 
$u\in L^2(\O)$ such that
\[
	\int_{\O}u(x)f(x)dx=-\int_{\G}g(x)\frac{\d v_f}{\d\nu_L}(x)\ds_x.
\]
  Moreover, by the above equation and \eqref{eq:dvf}, we have
\[
	|(u,f)|
\le 	\|g\|_{H^{-1/2}(\G)}\|\d_{\nu_L}v_f\|_{H^{1/2}(\G)}
\le	C\|g\|_{H^{-1/2}(\G)}\|f\|_{L^{2}(\O)},
\]
  for any $f\in L^2(\O)$, from which estimate \eqref{eq:ug'} follows.
\end{proof}

\bigskip
  Let $\Lm$ be a linear map which maps $g$ to the unique weak solution of 
\eqref{eq:ell}.
  Then we have seen that
\[
	\Lm\in\Lc(H^{3/2}(\G);H^2(\O))\cap\Lc(H^{-1/2}(\G);L^2(\O)),
\]
  where $\Lc(X;Y)$ denotes a set of bounded linear operators from a Banach 
space $X$ to another one $Y$.
  Then by the interpolation (see Theorem 5.1 in Chapter 1 of \cite{lm}), we 
have 
\[
	\Lm\in\Lc([H^{3/2}(\G),H^{-1/2}(\G)]_{\th};[H^2(\O);L^2(\O)]_{\th}).
\]
  for any $\th\in[0,1]$.
  In particular, we choose $\th=3/4$ and obtain
\[
	\Lm\in\Lc(L^2(\G);H^{1/2}(\O)).
\]
  That is, for any $g\in L^2(\G)$, $\Lm G$ belongs to $H^{1/2}(\O)$ and 
satisfies
\[
	\|\Lm g\|_{H^{1/2}(\O)}\le C\|g\|_{L^2(\G)}.
\]
  In particular, for any $0\le\th<1/4$, $\Lm g$ belongs to $\D(L^{\th})$ and 
satisfies
\begin{equation}\label{eq:Lg2}
	\|\Lm g\|_{\D(L^{\th})}\le C\|g\|_{H^{3/2}(\G)}.
\end{equation}
  By substituting $f=\v$ in \eqref{eq:weak}, we obtain
\begin{equation}\label{eq:Lgl}
	\l_n(\Lm g,\v)=-\<g,\d_{\nu_L}\v\>, \quad n=1,2,\dots.
\end{equation}

  For the arguments used here, we can refer to Chapter 2 in \cite{lm}, in 
which more general elliptic operator of order $2m$ is dealt with by assuming 
$C^{\infty}$-regularity for the coefficients $a_{ij}$ and the boundary $\G$.

  In Section \ref{sec:reg}, we apply the above results to the calculation of 
the eigenfunction expansion for the solution of \eqref{eq:sim}.

\section*{Acknowledgements}

  The author appreciates his supervisor Professor Masahiro Yamamoto for the 
useful advices.
  The author was granted by the Global COE program and is now supported by the 
FMSP program at Graduate School of Mathematical Sciences of The University of 
Tokyo.


\begin{thebibliography}{99}
	\bibitem{ag} E.E. Adams, L.W. Gelhar,
	Field study of dispersion in a heterogeneous aquifer 2. 
	Spatial moments analysis,
	Water Resources Res. 28 (1992) 3293-3307.

	\bibitem{ar} R.A. Adams,
	Sobolev Spaces, Academic Press,
	New York, 1975.

	\bibitem{aw} O.P. Agarwal,
	Solution for a fractional diffusion-wave equation defined in a bounded 
	domain,
	Nonlinear Dynam. 29 (2002) 145-155.

	\bibitem{am} S. Agmon, 
	Lectures on elliptic boundary value problems. Vol. 2.
	Am. Math. Soc., 1965.

	\bibitem{co} J.M. Coron,
	Control and Nonlinearity,
	Mathematical Surveys and Monographs 136. Am. Math. Soc., 
	Providence (2007)

	\bibitem{fy} K. Fujishiro, M. Yamamoto,
	Approximate controllability for fractional diffusion equations by 
	interior control,
	Appl. Anal. 93 (2014), no. 9, 1793-1810.

	\bibitem{fu} D. Fujiwara,
	Concrete characterization of the domains of fractional powers of some 
	elliptic differential operators of the second order, 
	Proceedings of the Japan Academy 43 (1967), no. 2, 82-86.

	\bibitem{gj} V.D. Gejji, H. Jafari,
	Boundary value problems for fractional diffusion-wave equation,
	Aust. J. Math. Anal. Appl. 3 (2006) 1-8.

	\bibitem{ha} Y. Hatano, N. Hatano,
	Dispersive transport of ions in column experiments: an explanation of 
	long-tailed profiles,
	Water Resources Res. 34 (1998) 1027-1033.

	\bibitem{is} V. Isakov,
	Inverse Problems for Partial Differential Equations,
	Springer-Verlag, Berlin, 2006.

	\bibitem{ks} A.A. Kilbas, H.M. Srivastava, J.J. Trujillo,
	Theory and Applications of Fractional Differential Equations,
	Elsevier, Amsterdam, 2006.

	\bibitem{li} J.L. Lions, 
	Remarks on approximate controllability, 
	J. Analyse Math., 59 (1992), 103-116.

	\bibitem{lm} J.L. Lions, E. Magenes, 
	Non-homogeneous Boundary Value Problems and Applications, vols. I, 
	Springer-Verlag, Berlin, 1972.

	\bibitem{lu} Y. Luchko,
	Some uniqueness and existence results for the initial-boundary value 
	problems for the generalized time-fractional diffusion equation,
	Comput. Math. Appl. 59 (2010) 1766-1772.

	\bibitem{ma} R.C. MacCamy, V.J. Mizel, T.I. Seidman,
	Approximate boundary controllability of the heat equation,
	J. Math. Anal. Appl., 23 (1968), 699-703.

	\bibitem{mk} R. Metzler, J. Klafter,
	The random walk's guide to anomalous diffusion: 
	a fractional dynamics approach, Physics reports 339.1 (2000): 1-77.

	\bibitem{mz} S. Micu, E. Zuazua, 
	An introduction to the controllability of partial differential 
	equations, 
	Quelques questions de th\'eorie du contr\^ole. 
	In: Sari T.(ed.) Collection Travaux en Cours (2004): 69-157.

	\bibitem{pd} I. Podlubny,
	Fractional Differential Equations,
	Academic Press, San Diego, 1999.

	\bibitem{ru} D.L. Russell,
	Controllability and stabilizability theory for linear partial 
	differential equations: recent progress and open questions,
	SIAM Rev. 20, (1978), pp. 639-739

	\bibitem{sk} K. Sakamoto, M. Yamamoto,
	Initial value/boundary value problems for fractional diffusion-wave 
	equations and applications to some inverse problems,
	J. Math. Anal. Appl. (2011).

	\bibitem{sm} S.G. Samko, A.A. Kilbas, O.I. Marichev,
	Fractional Integrals and Derivatives,
	Gordon and Breach Science Publishers, Philadelphia, 1993.

	\bibitem{sw} E.J.P.G. Schmidt, N. Weck,
	On the boundary behavior of solutions to elliptic and parabolic 
	equations---with applications to boundary control for parabolic 
	equations,
	SIAM J. Control Optim. 16 (1978) 593-598.

	\bibitem{st} E.M. Stein,
	Singular Integrals and Differentiability Properties of Functions,
	Princeton University Press, Princeton, 1970.

	\bibitem{wa} D. Washburn,
	A bound on the boundary input map for parabolic equations with 
	application to time optimal control, 
	SIAM J. Control Optim. 17 (1979), no. 5, 652-671. 

	\bibitem{yg} A. Yagi, 
	$H_{\infty}$ Functional Calculus and Characterization of Domains of 
	Fractional Powers, 
	Operator Theory: Advances and Applications, vol. 187 (2008), 217-235.

	\bibitem{zu} E. Zuazua, 
	Controllability and observability of partial differential equations: 
	Some results and open problems, 
	in: Handbook of Differential Equations: Evolutionary Differential
	Equations, vol. 3, Elsevier Science, 2006, pp. 527-621.
\end{thebibliography}
\end{document}